\documentclass[10pt]{article}
\usepackage{latexsym}
\usepackage{amsfonts}
\usepackage{enumerate}
\usepackage{multicol}
\usepackage{graphicx}
\usepackage{amssymb}
\usepackage{amsmath}
\usepackage{epic}

\title{The $h$-critical number of finite abelian groups \\[.4in]}

\author{B\'{e}la Bajnok \\[.1in] Department of Mathematics, Gettysburg College \\
Gettysburg, PA 17325-1486 USA \\E-mail:  bbajnok@gettysburg.edu \\[.4in]}

\date{December 4, 2014}

\newtheorem{thm}{Theorem}

\newtheorem{lem}[thm]{Lemma}
\newtheorem{cor}[thm]{Corollary}
\newtheorem{prop}[thm]{Proposition}
\newtheorem{conj}[thm]{Conjecture}

\begin{document}

\maketitle

\begin{abstract}
For a finite abelian group $G$ and a positive integer $h$, the unrestricted (resp.~restricted) $h$-critical number $\chi(G,h)$ (resp.~$\chi \hat{\;}(G,h)$) of $G$ is defined to be the minimum value of $m$, if exists, for which the $h$-fold unrestricted (resp.~restricted) sumset of every $m$-subset of $G$ equals $G$ itself.  Here we determine $\chi(G,h)$ for all $G$ and $h$; and prove several results for $\chi \hat{\;}(G,h)$, including the cases of any $G$ and $h = 2$, any $G$ and large $h$, and any $h$ for the cyclic group $\mathbb{Z}_n$ of even order.  We also provide a lower bound for $\chi \hat{\;}(\mathbb{Z}_n,3)$ that we believe is exact for every $n$---this conjecture is a generalization of the one made by Gallardo, Grekos, et al.~that was proved (for large $n$) by Lev.

\end{abstract}



\thispagestyle{empty}

\section{Introduction}

Throughout this paper, $G$ denotes a finite abelian group of order $n \geq 2$, written in additive notation.  For a positive integer $h$ and a nonempty subset $A$ of $G$, we let $hA$ and $h\hat{\;}A$ denote the {\em $h$-fold unrestricted sumset} and the {\em $h$-fold restricted sumset} of $A$, respectively; that is, $hA$ is the collection of sums of $h$ not-necessarily-distinct elements of $A$, and $h\hat{\;}A$ consists of all sums of $h$ distinct elements of $A$.  Furthermore,  we set $\Sigma A=\cup_{h=0}^{\infty} h \hat{\;} A$.    

The study of critical numbers originated with the 1964 paper \cite{ErdHei:1964a} of Erd\H{o}s and Heilbronn, in which they asked for the least integer $m$ so that for every set $A$ consisting of $m$ nonzero elements of the cyclic group $\mathbb{Z}_p$ of prime order $p$, we have $\Sigma A=\mathbb{Z}_p$.  More generally, one can define the {\em critical number} of $G$ as
$$\stackrel{*}{\chi} \hat{\;} (G)=\min \{m \; : \; A \subseteq G \setminus \{0\},  |A| \geq m \Rightarrow \Sigma  A=G \}.$$
Here the $\stackrel{*}{\;}$ indicates that only subsets of $G \setminus \{0\}$ are considered; alternately, some have studied 
$${\chi} \hat{\;} (G)=\min \{m \; : \; A \subseteq G,  |A| \geq m \Rightarrow \Sigma  A=G \}.$$

It took nearly half a century, but now, due to the combined results of Diderrich and Mann \cite{DidMan:1973a}, Diderrich \cite{Did:1975a},  Mann and Wou \cite{ManWou:1986a}, Dias Da Silva and Hamidoune \cite{DiaHam:1994a}, Gao and Hamidoune \cite{GaoHam:1999a}, Griggs \cite{Gri:2001a}, and  Freeze, Gao, and Geroldinger \cite{FreGaoGer:2009a, FreGaoGer:2014a}, we have the critical number of every group:

\begin{thm} [The combined results of authors above] \label{combined critical}

Suppose that $n \geq 10$, and let $p$ be the smallest prime divisor of $n$.  Then
$$\stackrel{*}{\chi} \hat{\;} (G)={\chi} \hat{\;} (G)-1=\left\{
\begin{array}{ll}
\lfloor 2 \sqrt{n-2} \rfloor & \mbox{if $G$ is cyclic of order $n=p$ or $n=pq$ where} \\
& \mbox{$q$ is prime and $3 \leq p \leq q \leq p+\lfloor 2 \sqrt{p-2} \rfloor+1$}\footnotemark, \\ \\ 
n/p+p-2 & \mbox{otherwise}.
\end{array}
\right.
$$

\end{thm}
\footnotetext{Note that $\lfloor 2 \sqrt{n-2} \rfloor=n/p+p-1$ in this case.}
We note that considering unrestricted sums rather than restricted sums makes the problem trivial: the corresponding unrestricted critical numbers ${\chi} (G)$ and $\stackrel{*}{\chi}(G)$, using the notations of Theorem \ref{combined critical}, are clearly given by
$$\stackrel{*}{\chi}(G)={\chi} (G)-1=n/p.$$

We now turn to our present subject: the critical number when only a fixed number of terms are added.  Here we consider both unrestricted sums and restricted sums; in particular, for a positive integer $h$, we define---if they exist, more on this below---the {\em unrestricted $h$-critical number} $\chi (G,h)$ and the {\em restricted $h$-critical number} $\chi \hat{\;} (G,h)$ as the minimum values of $m$ for which, respectively, the $h$-fold sumset and the $h$-fold restricted sumset of every $m$-element subset of $G$ is $G$ itself:
$${\chi}  (G, h)=\min \{m \; : \; A \subseteq G,  |A| \geq m \Rightarrow h  A=G \},$$
$${\chi} \hat{\;} (G, h)=\min \{m \; : \; A \subseteq G,  |A| \geq m \Rightarrow h \hat{\;} A=G \}.$$  

It is easy to see that for all $G$ and $h$ we have $hG=G$, so ${\chi}  (G, h)$ is always well defined; in Section \ref{prelim} we determine that ${\chi} \hat{\;} (G, h)$ is  well defined if, and only if, one of the following holds:
\begin{itemize}
  \item $h \in \{1,n-1\}$,
  \item $h \in \{2,n-2\}$, and $G$ is not isomorphic to an elementary abelian 2-group,
  \item $3 \leq h \leq n-3$. 
\end{itemize}
Furthermore, in Section \ref{prelim} we explain that the versions 
$$\stackrel{*}{\chi}  (G, h)=\min \{m \; : \; A \subseteq G \setminus \{0\},  |A| \geq m \Rightarrow h  A=G \}$$ and
$$\stackrel{*}{\chi} \hat{\;} (G, h)=\min \{m \; : \; A \subseteq G \setminus \{0\},  |A| \geq m \Rightarrow h \hat{\;} A=G \}$$
need not be studied separately, since---other than some trivial cases that we specify---they are well defined whenever their non-$\stackrel{*}{\;}$ versions are, and we have
$$\stackrel{*}{\chi}  (G, h)={\chi}  (G, h)$$ and $$\stackrel{*}{\chi} \hat{\;} (G, h)={\chi} \hat{\;} (G, h).$$

So let us see what we can say about the quantities ${\chi}  (G, h)$ and ${\chi} \hat{\;} (G, h)$.  We can determine the exact value of ${\chi}  (G, h)$, as follows.  

Recall that the minimum size $$\rho(G, m, h) = \min \{ |hA| \; : \; A \subseteq G, |A|=m\}$$ of $h$-fold sumsets of $m$-subsets of $G$ is known for all $G$, $m$, and $h$.  To state the result, we need the function
$$u(n,m,h)=\min \{f_d \; : \; d \in D(n)\},$$ where $n$, $m$, and $h$ are positive integers, $D(n)$ is the set of positive divisors of $n$, and
$$f_d=\left(h\left \lceil m/d \right \rceil-h +1 \right) \cdot d.$$  (Here $u(n,m,h)$ is a relative of the Hopf--Stiefel function used also in topology and bilinear algebra; see, for example, \cite{Sha:1984a}, \cite{Pla:2003a}, and \cite{Kar:2006a}.)  We then have:  

\begin{thm} [Plagne; cf.~\cite{Pla:2006a}] \label{value of u}
Let $n$, $m$, and $h$ be positive integers with $m \leq n$.  For any abelian group $G$ of order $n$ we have
$$\rho (G, m, h)=u(n,m,h).$$
\end{thm}   

Theorem \ref{value of u} allows us to determine ${\chi}  (G, h)$; in order to do so, we introduce a---perhaps already familiar---function first.   

Suppose that $h$ and $g$ are fixed positive integers; since we will only need the cases when $1 \leq g \leq h$, we make that assumption here.  Recall that we let $D(n)$ denote the set of positive divisors of $n$.  We then define  \label{vdef}
$$v_g(n,h)= \max \left\{ \left( \left \lfloor \frac{d-1-\mathrm{gcd} (d, g)}{h} \right \rfloor +1  \right) \cdot \frac{n}{d}  \; : \; d \in D(n) \right\}.$$  

We should note that the function $v_g(n,h)$ has appeared elsewhere in additive combinatorics already.  For example, according to the classical result of 
Diamanda and Yap (see \cite{DiaYap:1969a}), the maximum size of a sum-free set (that is, a set $A$ that is disjoint from $2A$) in the cyclic group $\mathbb{Z}_n$ is given by
$$v_1(n,3)= \left\{
\begin{array}{ll}
\left(1+\frac{1}{p}\right) \frac{n}{3} & \mbox{if $n$ has prime divisors congruent to 2 mod 3,} \\ & \mbox{and $p$ is the smallest such divisor,}\\ \\
\left\lfloor \frac{n}{3} \right\rfloor & \mbox{otherwise;}\\
\end{array}\right.$$ similarly, this author proved (see \cite{Baj:2009a}) that the maximum size of a $(3,1)$-sum-free set in $\mathbb{Z}_n$ (where $A$ is disjoint from $3A$) 
equals 
$$v_2(n,4) =\left\{
\begin{array}{ll}
\left(1+\frac{1}{p}\right) \frac{n}{4} & \mbox{if $n$ has prime divisors congruent to 3 mod 4,} \\ & \mbox{and $p$ is the smallest such divisor,}\\ \\
\left\lfloor \frac{n}{4} \right\rfloor & \mbox{otherwise.}\\
\end{array}\right.$$ It is believed that the analogous result for $(k,l)$-sum-free sets in $\mathbb{Z}_n$ (where $kA \cap lA=\emptyset$ for positive integers $k>l$) is given by $v_{k-l}(n,k+l)$; this was established for the case when $k-l$ and $n$ are relatively prime by Hamidoune and Plagne (see \cite{HamPla:2003a}).  
In Section \ref{nu formula} we provide a simpler alternate formula for $v_g(n,h)$, from which the expressions for $v_1(n,3)$ and $v_2(n,4)$ above will readily follow.

Returning now to the $h$-critical number of groups, in Section \ref{critU} we prove that for every group $G$ of order $n$ and for every $h$, we have
$${\chi}  (G, h)=v_1(n,h)+1.$$

Evaluating the restricted $h$-critical number ${\chi} \hat{\;} (G, h)$ seems much more challenging, and this is, of course, due to the fact that we do not have a general formula for 
the minimum size $$\rho \hat{\;} (G, m, h) = \min \{ |h \hat{\;}  A| \; : \; A \subseteq G, |A|=m\}$$ of $h$-fold restricted sumsets of $m$-subsets of $G$.  Indeed, we do not even know the value of $\rho \hat{\;} (G, m, h)$ for cyclic groups $G$ and $h=2$.  Essentially the only general result is for groups of prime order; solving a conjecture made by Erd\H{o}s and Heilbronn three decades earlier---not mentioned in \cite{ErdHei:1964a} but in \cite{ErdGra:1980a}---Dias Da Silva and Hamidoune succeeded in proving the following:

\begin{thm} [Dias Da Silva and Hamidoune; cf.~\cite{DiaHam:1994a}] \label{Dias Da Silva and Hamidoune}
For a prime $p$ and integers $1 \leq h \leq m \leq p$,  we have $$\rho\hat{\;} (\mathbb{Z}_p,m,h) = \mathrm{min} \{p, hm-h^2+1\}.$$

\end{thm}  
(The result was reestablished, using different methods, by Alon, Nathanson, and Ruzsa; see \cite{AloNatRuz:1995a}, \cite{AloNatRuz:1996a}, and \cite{Nat:1996a}.)  As a consequence, we have:

\begin{cor}  \label{prime rest crit}
For any positive integer $h$ and prime $p$ with $h \leq p-1$ we have
$$\chi \hat{\;} (\mathbb{Z}_p,h)=\left \lfloor (p-2)/h \right \rfloor +h+1.$$

\end{cor}

Let us see what else we can say about $\chi \hat{\;} (G,h)$.  Trivially, for all groups $G$ of order $n$ we have
$$\chi \hat{\;} (G,1)=\chi \hat{\;} (G,n-1)=n.$$
In Section \ref{critR h=2 and big}, we prove that for all $G$ of order $n$ and exponent at least 3, we have
$$\chi \hat{\;} (G,2)= (n+|L|)/2+1,$$ where $L$ denotes the subgroup of $G$ that consists of elements of order at most 2.  (Note that $n+|L|$ is always even; note also that for a group of exponent 2, $n=|L|$.)   In particular, for $n \geq 3$ we have $$\chi \hat{\;} (\mathbb{Z}_n,2) = \lfloor n/2 \rfloor +2.$$ 
As a consequence, we also show that this implies that if $G$ has order $n$ and exponent at least 3, and $h$ is an integer with $$(n+|L|)/2-1 \leq h \leq n-2,$$
then $$\chi \hat{\;} (G,h) = h+2.$$

This leaves us with the task of determining $\chi \hat{\;} (G,h)$ for groups of composite order and $$3 \leq h \leq (n+|L|)/2-2.$$  
In Section \ref{critR even n} we complete this task for cyclic groups of even order; namely, we prove that for an even value of $n \geq 12$, we have
$$\chi \hat{\;} (\mathbb{Z}_{n},h) = \left \{
\begin{array}{cl}
n/2+1 & \mbox{if} \; 3 \leq h \leq  n/2-2; \\ \\
n/2+2 & \mbox{if} \; h=n/2-1.
\end{array}
\right.$$
(This result was established for $h=3$ by Gallardo, Grekos, et al.~in ~\cite{GalGre:2002a}; our proof for the general case is based on their method.)

In Section \ref{case of h=3} we take a closer look at the case of $h=3$.  First, we prove tight lower bounds for $\chi \hat{\;} (\mathbb{Z}_n,3)$ as $n \geq 11$.  Namely, if $n$ has prime divisors congruent to $2$ mod $3$ and $p$ is the smallest such divisor, then we show that
  $$\chi \hat{\;} (\mathbb{Z}_n,3) \geq 
\left\{
\begin{array}{ll}
\left(1+\frac{1}{p} \right) \frac{n}{3} +3 & \mbox{if} \; n=p \; \mbox{or} \; n=15,  \\ \\
\left(1+\frac{1}{p} \right) \frac{n}{3} +2 & \mbox{if} \; n=3p \; \mbox{with} \; p \neq 5,  \\ \\
\left(1+\frac{1}{p} \right) \frac{n}{3} +1 & \mbox{otherwise}; 
\end{array} \right.$$
and if $n$ has no prime divisors congruent to $2$ mod $3$, then we prove that
  $$\chi \hat{\;} (\mathbb{Z}_n,3) \geq 
\left\{
\begin{array}{ll}
\left \lfloor \frac{n}{3} \right \rfloor +4 & \mbox{if $n$ is divisible by $9$},  \\ \\
\left \lfloor \frac{n}{3} \right \rfloor +3 & \mbox{otherwise}. 
\end{array} \right.$$
We also claim that, actually, equality holds above for all $n$---this is certainly the case if $n$ is even or prime; we have verified this (by computer) for all $n \leq 50$; and in Section \ref{case of h=3} we prove that equality follows from a conjecture that appeared in \cite{Baj:2013a}.  Our conjecture is a generalization of the one made by Gallardo, Grekos, et al.~in \cite{GalGre:2002a} that was proved (for large $n$) by Lev in \cite{Lev:2002a}.  

The pursuit of finding the value of $\chi \hat{\;} (G,h)$ in general remains challenging and exciting.

\section{Preliminary results} \label{prelim}

In this section we establish the conditions under which the four quantities ${\chi}  (G, h)$, ${\chi} \hat{\;} (G, h)$, $\stackrel{*}{\chi}  (G, h)$, and $\stackrel{*}{\chi} \hat{\;} (G, h)$ exist; furthermore, we show that, when they exist, then 
$${\chi}  (G, h)=\stackrel{*}{\chi}  (G, h)$$ and
$${\chi} \hat{\;} (G, h)=\stackrel{*}{\chi} \hat{\;} (G, h).$$

We start with the following easy result: 

\begin{prop}  \label{rho(G,m,h)>=m}

Let $A$ be an $m$-subset of $G$ and $h$ be a positive integer.  

\begin{enumerate} 
  \item If either \begin{enumerate} \item $h=1$ or 
  \item $A$ is a coset of a subgroup of $G$,    
  \end{enumerate} then $|h  A| =m.$
\item In all other cases, $|h  A| \geq m+1.$
\end{enumerate}  
\end{prop}

{\em Proof}:  The first claim is trivial.   To prove the second claim, we assume that $h \geq 2$ and that $|hA| \leq |A|=m$.  We will show that for any $a \in A$, we have $A=a+H$, where $H$ is the stabilizer subgroup of $(h-1)A$; that is,
$$H=\{g \in G  \mid g+(h-1)A=(h-1)A\}.$$

Consider the set $A'=A-a$.  Then $|A'|=m$ and $0 \in A'$, and therefore
$$(h-1)A =\{0\}+(h-1)A \subseteq A'+(h-1)A.$$
But then
$$|hA|=|hA-a|=|A'+(h-1)A|\geq |(h-1)A| \geq |(h-2) \cdot a+A|=|A|;$$
since we assumed $|hA| \leq |A|$, equality must hold throughout, and thus $$A'+(h-1)A = (h-1)A.$$ Therefore, $A' \subseteq H$, and so $A \subseteq a+H$, which implies that
$$|a+H| \geq |A| \geq |hA| =|(h-1)A+A| \geq |(h-1)A|=|H+(h-1)A| \geq |H|=|a+H|.$$  Then equality must hold throughout, and thus $a+H=A$, establishing our claim.
$\Box$  

As an immediate corollary, we see that ${\chi}  (G, h)$ is well defined for all $G$ and $h$, and $\stackrel{*}{\chi}  (G, h)$ is well defined if, and only if, the trivial conditions $n \geq 3$ and $h \geq 2$ hold.  

The version of Proposition \ref{rho(G,m,h)>=m} for restricted sumsets is substantially more complicated:

\begin{thm} [Girard, Griffiths, and Hamidoune;  cf.~\cite{GirGriHam:2012a}] \label{rhohat(G,m,h)>=m}

Let $A$ be an $m$-subset of $G$, and suppose that $1 \leq h \leq m-1$.  We let $L$ denote the subgroup of $G$ that consists of elements of order at most 2.

\begin{enumerate} 
  \item If $h \in \{2,m-2\}$ and $A$ is a coset of a subgroup of $L$, then $|h  \hat{\;} A|=m-1.$
  \item If any of the conditions \begin{enumerate} \item $h \in \{1,m-1\}$, 
  \item $A$ is a coset of a subgroup of $G$,    
  \item $h \in \{2,m-2\}$ and $A$ consists of all but one element of a coset of a subgroup of $L$, or 
  \item $h \in \{2,m-2\}$ and $m=4$ and $A$ consists of two cosets of a subgroup of order 2
  \end{enumerate} holds, then $|h \hat{\;} A| =m.$
\item In all other cases, $|h \hat{\;} A| \geq m+1.$
\end{enumerate}  
\end{thm}

As a consequence, we get that ${\chi} \hat{\;} (G, h)$ is  well defined if, and only if, one of the following holds:
\begin{itemize}
  \item $h \in \{1,n-1\}$,
  \item $h \in \{2,n-2\}$, and $G$ is not isomorphic to an elementary abelian 2-group,
  \item $3 \leq h \leq n-3$; 
\end{itemize}
and $\stackrel{*} {\chi} \hat{\;} (G, h)$ is  well defined if, and only if, one of the following holds:
\begin{itemize}
  \item $n=5$ and $h=2$,
  \item $n \geq 6$, $h \in \{2,n-2\}$, and $G$ is not isomorphic to an elementary abelian 2-group;
  \item $3 \leq h \leq n-3$.
\end{itemize}
From this we can conclude that, other than the trivial cases of $h \in \{1,n-1\}$ or $n \leq 5$, $\stackrel{*} {\chi} \hat{\;} (G, h)$ is  well defined exactly when $ {\chi} \hat{\;} (G, h)$ is.

Next we prove that our $\stackrel{*}{\;}$ quantities are equal to their respective non-$\stackrel{*}{\;}$ versions:

\begin{prop}

When they are defined, we have
$$\stackrel{*}{\chi}  (G, h)={\chi}  (G, h)$$ and $$\stackrel{*}{\chi} \hat{\;} (G, h)={\chi} \hat{\;} (G, h).$$

\end{prop}

{\em Proof}: We only prove the first claim as the other is similar.  For that, the other direction being obvious, we just need to show that  
$$\stackrel{*}{\chi}  (G,h) \geq {\chi}  (G,h).$$
To see this, let $B$ be a subset of $G$ of size ${\chi}  (G,h)-1$ for which $h  B \neq G$.  Since $|B| \leq n-1$, we have $|-B| \leq n-1$ as well; let $g \in G \setminus (-B)$.  Then $A=g+B$ has size ${\chi}  (G,h)-1$, and $A \subseteq G \setminus \{0\}$, since $0 \in A$ would contradict $g \not \in -B$.  But $h  A$ and $h  B$ have the same size, so we conclude that $h  A \neq G$, from which 
our inequality follows. $\Box$

To summarize this section: it suffices to study ${\chi} (G, h)$ and ${\chi} \hat{\;} (G, h).$

\section{The function $v_g(n,h)$} \label{nu formula}

In this section we prove a result that greatly simplifies the evaluation of the function $v_g(n,h)$ that we defined in the Introduction.

\begin{thm} \label{nu}
Suppose that $n$, $h$, and $g$ are positive integers and that $1 \leq g \leq h$.  For $i=2,3, \dots,h-1$, let $P_i(n)$ be the set of those prime divisors of $n$ that do not divide $g$ and that leave a remainder of $i$ when divided by $h$; that is, 
$$P_i(n)=\{ \; p \in D(n) \setminus D(g)  \; : \; p \; \mbox{prime and} \; p \equiv i \; (\mathrm{mod} \; h) \}.$$  We let $I$ denote those values of $i=2,3,\dots,h-1$ for which $P_i(n) \not = \emptyset$, and for each $i \in I$, we let $p_i$ be the smallest element of $P_i(n)$.

Then, the value of $v_g(n,h)$ is
$$v_g(n,h) =\left\{
\begin{array}{lll}
\frac{n}{h}  \cdot \mathrm{max} \left\{ 1+ \frac{h-i}{p_i}  \; : \; i \in I \right\} & \mbox{if} & I \not = \emptyset; \\ \\
\left\lfloor \frac{n}{h} \right\rfloor & \mbox{if} & I = \emptyset \; \mbox{and} \; g \not = h;\\
& & \\ 
\left\lfloor \frac{n-1}{h} \right\rfloor & \mbox{if} & I = \emptyset \; \mbox{and} \; g = h. \\
\end{array}\right.$$
 
\end{thm} 

{\em Proof}:  Suppose that $d$ is a positive divisor of $n$, and    define the function
$$f(d)= \left(  \left \lfloor \frac{d-1-\mathrm{gcd}(d,g)}{h} \right \rfloor +1 \right) \cdot \frac{n}{d}.$$

We first prove the following.

{\bf Claim 1:}  Let $i$ be the remainder of $d$ when divided by $h$.  We then have
$$ f(d) =
\left\{
\begin{array}{ll}
\frac{n}{h} \cdot \left(1+ \frac{h-i}{d} \right) & \mbox{if } \; \mathrm{gcd}(d,g) < i; \\ \\
\frac{n}{h}\cdot \left( 1-\frac{h}{d} \right) & \mbox{if } \; h|d \; \mbox{and} \; g=h; \\ \\
\frac{n}{h} \cdot \left( 1-\frac{i}{d} \right) & \mbox{otherwise. }
\end{array}\right.$$

{\em Proof of Claim 1.}  We start with
$$ \left \lfloor \frac{d-1-\mathrm{gcd}(d,g)}{h} \right \rfloor=\frac{d-i}{h}+\left \lfloor \frac{i-1-\mathrm{gcd}(d,g)}{h} \right \rfloor.$$
We investigate the maximum and minimum values of the quantity $\left \lfloor \frac{i-1-\mathrm{gcd}(d,g)}{h} \right \rfloor$.

For the maximum, we have 
$$\left \lfloor \frac{i-1-\mathrm{gcd}(d,g)}{h} \right \rfloor \leq \left \lfloor \frac{(h-1)-1-1}{h} \right \rfloor \leq 0,$$ with equality if, and only if, $i-1-\mathrm{gcd}(d,g) \geq 0$; that is, $\mathrm{gcd}(d,g) < i$.

For the minimum, we get
$$\left \lfloor \frac{i-1-\mathrm{gcd}(d,g)}{h} \right \rfloor \geq \left \lfloor \frac{0-1-g}{h} \right \rfloor \geq 
\left \lfloor \frac{0-1-h}{h} \right \rfloor =-2,$$
with equality if, and only if, $i=0$, $\mathrm{gcd}(d,g)=g$, and $g=h$; that is, $h|d$ and $g=h$.  

The proof of Claim 1 now follows easily.  $\Box$

{\bf Claim 2:}  Using the notations as above, assume that $\mathrm{gcd}(d,g) \geq i$.  Then
$$ f(d) \leq
\left\{
\begin{array}{ll}
\frac{n}{h}  & \mbox{if } \;  g \neq h; \\ \\
\frac{n-1}{h}& \mbox{if } \;  g = h.
\end{array}\right.$$

{\em Proof of Claim 2.}  By Claim 1, we have 
$$f(d) \leq \frac{n}{h}.$$    
Furthermore, unless $i=0$ and $g \neq h$, we have 
$$f(d) \leq \frac{n}{h}\cdot \left( 1-\frac{1}{d} \right) \leq \frac{n}{h}\cdot \left( 1-\frac{1}{n} \right) = \frac{n-1}{h}.$$
$\Box$

{\bf Claim 3:}  For all $g$, $h$, and $n$ we have 
$$ v_g(n,h) \geq
\left\{
\begin{array}{ll}
\left \lfloor \frac{n}{h} \right \rfloor & \mbox{if } \; g \neq h; \\ \\
\left \lfloor \frac{n-1}{h} \right \rfloor & \mbox{if } \; g=h.
\end{array}\right.$$

{\em Proof of Claim 3.}  We first note that 
\begin{eqnarray*}
v_g(n,h) &=& \max \left\{ \left( \left \lfloor \frac{d-1-\mathrm{gcd} (d, g)}{h} \right \rfloor +1  \right) \cdot \frac{n}{d}  \; : \; d \in D(n) \right\} \\ \\
& \geq & \left \lfloor \frac{n-1-\mathrm{gcd} (n, g)}{h} \right \rfloor +1 \\ \\
& \geq & \left \lfloor \frac{n-1-g}{h} \right \rfloor +1.
\end{eqnarray*}
The claim now follows, since $g+1 \leq h$, unless $g=h$ in which case $$\left \lfloor \frac{n-1-g}{h} \right \rfloor +1=\left \lfloor \frac{n-1}{h} \right \rfloor.$$   $\Box$

We are now ready for the proof of Theorem \ref{nu}.

{\em Proof of Theorem \ref{nu}.}  Let $d_0$ be any positive divisor of $n$ for which $ v_g(n,h) =f(d_0)$; let $i_0$ be the remainder of $d_0$ mod $h$.  The following two claims together establish Theorem \ref{nu}.  

{\bf Claim 4:}  If $\mathrm{gcd}(d_0,g) \geq i_0$, then $I=\emptyset$ and $$ v_g(n,h) =
\left\{
\begin{array}{ll}
\left \lfloor \frac{n}{h} \right \rfloor & \mbox{if } \; g \neq h; \\ \\
\left \lfloor \frac{n-1}{h} \right \rfloor & \mbox{if } \; g=h
\end{array}\right.$$

{\em Proof of Claim 4:}  By Claim 2, $$ v_g(n,h) =f(d_0) \leq n/h.$$  If we were to have an element $i \in I$, then for the corresponding prime divisor $p_i$ of $n$ we have $$\mathrm{gcd}(p_i,g)=1 < i,$$ thus by Claim 1, 
$$v_g(n,h) \geq f(p_i) = \frac{n}{h} \cdot \left( 1+ \frac{h-i}{p_i} \right) > \frac{n}{h},$$ a contradiction.  The result now follows from Claims 2 and 3.  $\Box$

{\bf Claim 5:}  If $\mathrm{gcd}(d_0,g) < i_0$, then $i_0 \in I$, $d_0 \in P_{i_0} (n)$, and $$v_g(n,h) = \frac{n}{h} \cdot \left( 1+ \frac{h-i_0}{d_0} \right).$$

{\em Proof of Claim 5:}  First, we prove that $d_0$ is prime.  Note that our assumption implies that $i_0 \geq 2$, and thus $d_0$ has no divisor that is divisible by $h$, and has at least one prime divisor that leaves a remainder greater than 1 mod $h$.  Let $p$ be the smallest prime divisor of $d_0$ that leaves a remainder more than 1 mod $h$, and let $i$ be this remainder.

We establish the inequality
$$\frac{h-2}{p^2} < \frac{h-i}{p},$$ as follows.
Since $i \leq h-1$, the inequality clearly holds when $p>h-2$, so let us assume that $p\leq h-2$.  Note that, in this case, $i=p$, so we need to establish that $$\frac{h-2}{p^2} < \frac{h-p}{p};$$ this is not hard either since we have
$$h-2 =hp-h(p-1)-2 \leq hp-(p+2)(p-1)-2=hp-p^2-p<(h-p)p.$$

Assume now that $i \neq i_0$, and thus $d_0/p \not \equiv 1$ mod $h$.  Then $d_0/p$ also has a prime divisor, say $p'$, that leaves a remainder greater than 1 mod $h$, and by the choice of $p$, $p' \geq p$ and thus $d_0 \geq p^2$.  But then we have
$$v_g(n,h) =f(d_0)=\frac{n}{h} \cdot \left( 1+ \frac{h-i_0}{d_0} \right) \leq \frac{n}{h} \cdot \left( 1+ \frac{h-2}{p^2} \right)<\frac{n}{h} \cdot \left( 1+ \frac{h-i}{p} \right)=f(p),$$ a contradiction.
 
Therefore, $i=i_0$, and thus 
$$v_g(n,h) =f(d_0)=\frac{n}{h} \cdot \left( 1+ \frac{h-i_0}{d_0} \right) \leq \frac{n}{h} \cdot \left( 1+ \frac{h-i_0}{p} \right)=f(p);$$ since we must have equality, $d_0=p$ follows.  

This establishes the fact that $d_0$ is prime.  Since $$\mathrm{gcd}(d_0,g) < i_0 \leq d_0,$$ $d_0$ cannot divide $g$.  This establishes Claim 5, and thus completes the proof of Theorem \ref{nu}.
$\Box$ 

We should also note that it is easy to show that, when $I \neq \emptyset$ in the statement of Theorem \ref{nu}, there is a unique $i$ (and thus $p_i$) for which $\frac{h-i}{p_i}$ is maximal.

\section{The unrestricted $h$-critical number} \label{critU}

Recall from our Introduction that 
$$v_1(n,h)= \max \left\{ \left( \left \lfloor \frac{d-2}{h} \right \rfloor +1  \right) \cdot \frac{n}{d}   \; : \; d \in D(n) \right\}.$$

Here we prove the following:

\begin{thm}  \label{h crit numb}
For all finite abelian groups $G$ of order $n$ and all positive integers $h$, the (unrestricted) $h$-critical number of $G$ equals $$\chi (G,h)=v_1(n,h)+1.$$
\end{thm}
   
{\em Proof}:  We need to prove that, for $m=v_1(n,h)$, we have $$u(n,m,h)<n$$ but $$u(n,m+1,h) \geq n.$$  Let $d_0 \in D(n)$ be such that $$v_1(n,h)=\max \left\{ \left( \left \lfloor \frac{d-2}{h} \right \rfloor +1  \right) \cdot \frac{n}{d}   \; : \; d \in D(n) \right\}=\left( \left \lfloor \frac{d_0-2}{h} \right \rfloor +1  \right) \cdot \frac{n}{d_0}.$$    

To establish the first inequality, simply note that $u(n,m,h) \leq f_{n/d_0}(m,h)$ where
$$f_{n/d_0}(m,h)=\left(h \cdot \left ( \left \lfloor \frac{d_0-2}{h} \right \rfloor +1 \right )    -h+1  \right) \cdot \frac{n}{d_0}=
\left(h \cdot \left \lfloor \frac{d_0-2}{h} \right \rfloor +1  \right) \cdot \frac{n}{d_0} \leq (d_0-1) \cdot \frac{n}{d_0} <n. $$

For the second inequality, we must prove that, for any $d \in D(n)$, we have $f_d(m+1,h) \geq n$; that is, 
$$h \cdot \left \lceil \frac{ 
\left( \left \lfloor \frac{d_0-2}{h} \right \rfloor +1 \right )  \cdot \frac{n}{d_0} +1}{d} 
\right \rceil  -h+1 \geq \frac{n}{d}.$$
 But $n/d \in D(n)$, so by the choice of $d_0$, we have
$$\left( \left \lfloor \frac{d_0-2}{h} \right \rfloor +1 \right )  \cdot \frac{n}{d_0} \geq \left( \left \lfloor \frac{n/d-2}{h} \right \rfloor +1 \right )  \cdot \frac{n}{n/d},$$ and thus
\begin{eqnarray*}
h \cdot \left \lceil \frac{ 
\left( \left \lfloor \frac{d_0-2}{h} \right \rfloor +1 \right )  \cdot \frac{n}{d_0} +1}{d} 
\right \rceil  -h+1
& \geq & h \cdot \left \lceil 
\left( \left \lfloor \frac{n/d-2}{h} \right \rfloor +1 \right )  +\frac{1}{d} 
\right \rceil  -h+1 \\
& = & h \cdot  
\left( \left \lfloor \frac{n/d-2}{h} \right \rfloor +2 \right )  
  -h+1 \\
& \geq & h \cdot  
\left(  \frac{n/d-2-(h-1)}{h}  +2 \right )  
  -h+1 \\
& = & \frac{n}{d}.
\end{eqnarray*}
Our proof is complete.  $\Box$

\section{The restricted $h$-critical number for $h=2$ and large $h$} \label{critR h=2 and big}

First, we evaluate $\chi \hat{\;} (G,2)$:

\begin{prop} \label{2-chrom G} 
Suppose that $G$ is of order $n$ and is not isomorphic to the elementary abelian 2-group, and let $L$ denote its subset---indeed, subgroup---consisting of elements of order at most 2.  Then $$\chi \hat{\;} (G,2) = (n+|L|)/2+1.$$  In particular, for $n \geq 3$ we have $$\chi \hat{\;} (\mathbb{Z}_n,2) = \lfloor n/2 \rfloor +2.$$

\end{prop}

We first prove the following.

\begin{lem} \label{lemma for doubling}
For a given $g \in G$, let $L_g=\{x \in G \mid 2x=g\}$.  If $L_g \neq \emptyset$, then $|L_g| = |L|.$
\end{lem}

{\em Proof}: Choose an element $x \in L_g$.  Then $x-L_g \subseteq L$, so $|x-L_g|=|L_g| \leq |L|.$  

Similarly, $x+L \subseteq L_g$, so  $|x+L|=|L| \leq |L_g|.$ $\Box$

{\em Proof of Proposition \ref{2-chrom G}}: Suppose first that $$m = (n+|L|)/2+1.$$ Note that our assumption on $G$ implies that $3 \leq m \leq n$.

Let $A$ be an $m$-subset of $G$, let $g \in G$ be arbitrary, and set $B= g-A$.  Then $|B|=m$, and thus
$$|A \cap B|=|A|+|B|-|A \cup B|  \geq 2m-n  = |L|+2.$$  
By our lemma above, we must have an element $a_1 \in A \cap B$ for which $a_1 \not \in L_g$.  Since $a_1 \in A \cap B$, we also have an element $a_2 \in A$ for which $a_1=g-a_2$ and thus $g=a_1+a_2$.  But $a_1 \not \in L_g$, and therefore $a_2 \neq a_1$.  In other words, $g \in 2\hat{\;}A$; since $g$ was arbitrary, we have $G=2\hat{\;}A$, as claimed.

For the other direction, we need to find a subset $A$ of $G$ with $$|A|=(n+|L|)/2$$ for which $2 \hat{\;} A \neq G$.  Observe that the elements of $G \setminus L$ are distinct from their inverses, so we have a (possibly empty) subset $K$ of $G \setminus L$ with which $$G=L \cup K \cup (-K),$$ and $L$, $K$, and $-K$ are pairwise disjoint.   Now set $A=L \cup K.$  Clearly, $A$ has the right size; furthermore, it is easy to verify that $0 \not \in 2\hat{\;}A$ and thus $2\hat{\;}A \neq G$.
$\Box$

Next, we show how Proposition \ref{2-chrom G} allows us evaluate $\chi \hat{\;} (G,h)$ for all large values of $h$.  In particular, we have:

\begin{prop}
Suppose that $G$ is not an elementary abelian 2-group and $h$ is a positive integer with $$(n+|L|)/2-1 \leq h \leq n-2.$$  Then $$\chi \hat{\;} (G,h) = h+2.$$

\end{prop}

{\em Proof}: Assume first that $A$ is an $(h+1)$-subset of $G$.  Then $$|h \hat{\;} A| =h+1 \leq n-1,$$ so $\chi \hat{\;} (G,h)$ is at least $h+2$.  

Now let $A$ be an $(h+2)$-subset of $G$.  Then, by symmetry, $|h \hat{\;} A| =|2 \hat{\;} A|;$ since 
$$|A|=h+2 \geq  (n+|L|)/2+1,$$  
by Proposition \ref{2-chrom G} we have $h \hat{\;} A =G.$  This establishes our claim.  $\Box$

\section{The restricted $h$-critical number of cyclic groups of even order} \label{critR even n}

Here we determine the value of $\chi \hat{\;} (\mathbb{Z}_{n},h)$ for all values of $h$ when $n$ is even:

\begin{thm}  \label{n even h critical}
Suppose that $n$ is even and $n \geq 12.$  Then 
$$\chi \hat{\;} (\mathbb{Z}_{n},h) = \left \{
\begin{array}{cl}
n & \mbox{if} \; h=1; \\ \\
n/2+2 & \mbox{if} \; h=2; \\ \\
n/2+1 & \mbox{if} \; h=3,4,\dots, n/2-2; \\ \\
n/2+2 & \mbox{if} \; h=n/2-1; \\ \\
h+2 & \mbox{if} \; h=n/2, n/2+1, \dots, n-2; \\ \\
n & \mbox{if} \; h=n-1.
\end{array}
\right.$$

\end{thm}

Theorem \ref{n even h critical} was established for $h=3$ by Gallardo, Grekos, et al.~in ~\cite{GalGre:2002a}; our proof for the general case is based on their method as well as Theorem \ref{rhohat(G,m,h)>=m} above.

{\em Proof}:  The  cases of $h \leq 2$ or $h \geq n/2$ have been already addressed, leaving only $3 \leq h \leq n/2-1$.  In fact, as we now show, it suffices to treat the cases of $3 \leq h \leq n/4$: 

To conclude that we then have $\chi \hat{\;} (\mathbb{Z}_{n},h) = n/2+1$ for $$n/4 +1 \leq h \leq n/2-2$$ as well, note that, obviously, $\chi \hat{\;} (\mathbb{Z}_{n},h) \geq n/2+1$, and that if $A$ is a subset of $\mathbb{Z}_n$ of size $n/2+1$, then, since $$3 \leq n/2+1-h \leq n/4,$$ we have
$$|h\hat{\;} A|=|(n/2+1-h) \hat{\;} A|=n.$$

Similarly, with $\chi \hat{\;} (\mathbb{Z}_{n},2) = n/2+2$ and $\chi \hat{\;} (\mathbb{Z}_{n},3) = n/2+1$ we can settle the case of $h=n/2-1$: Choosing a subset $A$ of $\mathbb{Z}_n$ of size $n/2+1$ for which $|2 \hat{\;}  A |<n$ implies that we also have $$|(n/2-1) \hat{\;}  A |<n$$ and thus $\chi \hat{\;} (\mathbb{Z}_{n},n/2-1)$ is at least $n/2+2$; while for any $B \subset \mathbb{Z}_n$ of size $n/2+2$ we get 
$$|(n/2-1)\hat{\;} B|=|3 \hat{\;} B|=n.$$
Therefore, for the rest of the proof, we assume that $3 \leq h \leq n/4$.  

Since we clearly have $\chi \hat{\;} (\mathbb{Z}_{n},h) \geq n/2+1$, it suffices to prove the reverse inequality.  For that, let $A$ be a subset of $\mathbb{Z}_n$ of size $n/2+1$; we need to prove that $h\hat{\;} A = \mathbb{Z}_n$.  

Let $O$ and $E$ denote the set of odd and even elements of $\mathbb{Z}_n$, respectively, and let $A_O$ and $A_E$ be the set of odd and even elements of $A$, respectively.  Note that both $A_O$ and $A_E$ have size at most $n/2$ and thus neither can be empty.  We will consider four cases:

Assume first that $|A_O| \leq 2$.  Then $|A_E| \geq n/2-1$.  Observe that $3 \leq h \leq n/4$ and $n \geq 12$ imply that $$2 \leq h-1 < h \leq n/2-3,$$ and $n/2-1$ is not a divisor of $n$.  Therefore, by Theorem \ref{rhohat(G,m,h)>=m}, both $(h-1)\hat{\;} A_E$ and $h\hat{\;} A_E$ have size at least $n/2$. But, of course, both $(h-1)\hat{\;} A_E$ and $h\hat{\;} A_E$ are subsets of $E$, so 
$$(h-1)\hat{\;} A_E=h\hat{\;} A_E=E.$$  

Now let $a$ be any element of $A_O$; we then see that $$a+(h-1)\hat{\;} A_E=a+E=O.$$  Therefore, $$(a+(h-1)\hat{\;} A_E) \cup h\hat{\;} A_E = O \cup E=\mathbb{Z}_n;$$ since both $a+(h-1)\hat{\;} A_E$ and $ h\hat{\;} A_E$ are subsets of $h\hat{\;} A$, we get $h\hat{\;} A = \mathbb{Z}_n$.

Next, we assume that $|A_E| \leq 2$.  In this case, an argument similar to the one in the previous case yields that
$$(h-1)\hat{\;} A_O = \left \{
\begin{array}{cl}
O & \mbox{if $h$ is even}, \\
E & \mbox{if $h$ is odd}; \\  
\end{array}
\right.$$
and 
$$h\hat{\;} A_O = \left \{
\begin{array}{cl}
E & \mbox{if $h$ is even}, \\
O & \mbox{if $h$ is odd}. \\  
\end{array}
\right.$$
Let $a$ be any element of $A_E$; we get
$$(a+(h-1)\hat{\;} A_O) \cup h\hat{\;} A_O = \mathbb{Z}_n$$ regardless of whether $h$ is even or odd; therefore, $h\hat{\;} A = \mathbb{Z}_n$.

Before turning to the last two cases, we observe that, since $h \leq n/4$, we have $$|A|=n/2+1 \geq 2h+1,$$ and thus at least one of $A_O$ or $A_E$ must have size at least $h+1$.  

Consider the case when $|A_O| \geq 3$ and $|A_E| \geq h+1$.  Referring to Theorem \ref{rhohat(G,m,h)>=m} again, we deduce that $(h-2)\hat{\;} A_E$ and $(h-1)\hat{\;} A_E$ both have size at least $|A_E|$, and that $2\hat{\;} A_O$ is of size at least $|A_O|$.  

Now let $g_O$ be any element of $O$; we have
$$|g_O -A_O|+|(h-1)\hat{\;} A_E| \geq |A_O|+|A_E|=n/2+1.$$ But $g_O -A_O$ and $(h-1)\hat{\;} A_E$ are both subsets of $E$, so they cannot be disjoint; this then means that $g_O$ can be written as the sum of an element of $A_O$ and $h-1$ distinct elements of $A_E$, so $g_O \in h\hat{\;} A$.   

Similarly, for any element $g_E$  of $E$, we have  
$$|g_E -(h-2)\hat{\;} A_E|+|2\hat{\;} A_O| \geq |A_E|+|A_O|=n/2+1,$$ and thus $g_E$ can be written as the sum of $h-2$ distinct elements of $A_E$ and two distinct elements of $A_O$, so $g_E \in h\hat{\;} A$.

Combining the last two paragraphs yields $O \cup E \subseteq h\hat{\;} A$ and thus $h\hat{\;} A = \mathbb{Z}_n$.

For our fourth case, assume that $|A_E| \geq 3$ and $|A_O| \geq h+1$.  As above, we can conclude that $|(h-2)\hat{\;} A_O| \geq |A_O|$, $|(h-1)\hat{\;} A_O| \geq |A_O|$, and $|2\hat{\;} A_E| \geq |A_E|$.

Let $g$ be any element of $\mathbb{Z}_n$.  If $g$ and $h$ are of the same parity (both even or both odd), then we find that $g-(h-2)\hat{\;} A_O$ and $2\hat{\;} A_E$ are each subsets of $E$.  As above, we see that they cannot be disjoint, and thus $$g \in (h-2)\hat{\;} A_O + 2\hat{\;} A_E \subseteq h\hat{\;} A.$$ 

The subcase when $g$ is even and $h$ is odd is similar: this time we see that $g-(h-1)\hat{\;} A_O$ and $A_E$ are each subsets of $E$ and that they cannot be disjoint, so $$g \in (h-1)\hat{\;} A_O + A_E \subseteq h\hat{\;} A.$$ 

The final subcase, when $g$ is odd and $h$ is even, needs more work.  We first prove that there is at most one element $a \in A_O$ for which $A_O \setminus \{a \}$ is the coset of a subgroup of $\mathbb{Z}_n$.  Suppose, indirectly, that $a_1$ and $a_2$ are distinct elements of $A_O$ so that $A_O \setminus \{a_1 \}$ and $A_O \setminus \{a_2 \}$ are both cosets.  In this case, they must be cosets of the same subgroup since $\mathbb{Z}_n$ has only one subgroup of that size.  But $|A_O| \geq 3$, so $A_O \setminus \{a_1 \}$ and $A_O \setminus \{a_2 \}$ are not disjoint, which implies that they are actually equal, which is a contradiction since $a_1$ is an element of $A_O \setminus \{a_2 \}$ but not of $A_O \setminus \{a_1 \}$.

We also need to consider the special case when $|A_O|=5$; we can then see that there is at most one element $a \in A_O$ for which $A_O \setminus \{a \}$ is the union of two cosets of $\{0, n/2\}$.         

Hence we have an element $a_O \in A_O$ so that $A_O \setminus \{a_O \}$ is not the coset of a subgroup of $\mathbb{Z}_n$, and not the union of two cosets of the subgroup of size 2.  But then, by Theorem \ref{rhohat(G,m,h)>=m}, 
$$|(h-2)\hat{\;} (A_O \setminus \{a_O \})| \geq |A_O|.$$  Therefore, $$|(h-2)\hat{\;} (A_O \setminus \{a_O \})|+|g-a_O- A_E| \geq |A_O|+|A_E|=n/2+1;$$ since both $(h-2)\hat{\;} (A_O \setminus \{a_O \})$ and $g-a_O- A_E$ are subsets of $E$, this can only happen if they are not disjoint, which means that $$g \in (h-2)\hat{\;} (A_O \setminus \{a_O \}) + (a_O+ A_E) \subseteq h \hat{\;} A.$$  This completes our proof.
$\Box$

\section{The restricted 3-critical number of cyclic groups} \label{case of h=3}

In this section we summarize what we can say about the case of $h=3$ in the cyclic group of order $n$.  

We will rely on the following result:

\begin{thm} [B.; cf.~\cite{Baj:2013a}] \label{cor h=3}
For all positive integers $n$ and $m$ with $4 \leq m \leq n$ we have
$$\rho\hat{\;} (\mathbb{Z}_n,m,3) \leq \left\{
\begin{array}{ll}
\min\{u(n,m,3), 3m-3-\gcd(n,m-1)\} & \mbox{if} \; \gcd(n,m-1) \geq 8; \\ \\
\min\{u(n,m,3), 3m-10\} & \mbox{if} \; \gcd(n,m-1) = 7, \; \mbox{or} \\
&  \gcd(n,m-1) \leq 5, \; 3|n, \; \mbox{and} \; 3|m, \; \mbox{or} \\
&  \gcd(n,m-1) \leq 5, \; (3m-9)|n, \; \mbox{and} \; 5|(m-3); \\ \\
\min\{u(n,m,3), 3m-9\} & \mbox{if} \; \gcd(n,m-1) = 6, \; \mbox{or} \\
& m=6 \; \mbox{and} \; 10|n \; \mbox{but} \; 3 \not | n; \\ \\
\min\{u(n,m,3), 3m-8\} & \mbox{otherwise.} 
\end{array}
\right.$$

\end{thm}

Our result for $\chi \hat{\;} (\mathbb{Z}_n,3)$ is, as follows:

\begin{prop} \label{h=3 rest crit}
Let $n$ be an arbitrary integer with $n \geq 11$.  
\begin{enumerate}
  \item If $n$ has prime divisors congruent to $2$ mod $3$ and $p$ is the smallest such divisor, then
  $$\chi \hat{\;} (\mathbb{Z}_n,3) \geq 
\left\{
\begin{array}{ll}
\left(1+\frac{1}{p} \right) \frac{n}{3} +3 & \mbox{if} \; n=p \; \mbox{or} \; n=15;  \\ \\
\left(1+\frac{1}{p} \right) \frac{n}{3} +2 & \mbox{if} \; n=3p \; \mbox{with} \; p \neq 5;  \\ \\
\left(1+\frac{1}{p} \right) \frac{n}{3} +1 & \mbox{otherwise}. 
\end{array} \right.$$
  \item If $n$ has no prime divisors congruent to $2$ mod $3$, then
  $$\chi \hat{\;} (\mathbb{Z}_n,3) \geq 
\left\{
\begin{array}{ll}
\left \lfloor \frac{n}{3} \right \rfloor +4 & \mbox{if $n$ is divisible by $9$};  \\ \\
\left \lfloor \frac{n}{3} \right \rfloor +3 & \mbox{otherwise}. 
\end{array} \right.$$
\end{enumerate}

\end{prop}

{\em Proof}: Note that the case when $n$ is even follows from Theorem \ref{n even h critical}, since $$\left(1+\frac{1}{2} \right) \frac{n}{3} +1=\frac{n}{2}+1;$$  and the case when $n$ is prime follows from Theorem \ref{prime rest crit} since
$$\left \lfloor \frac{p-2}{3} \right \rfloor +3+1 =
\left\{
\begin{array}{ll}
\left(1+\frac{1}{p} \right) \frac{p}{3} +3 & \mbox{if} \; p \equiv 2 \; \mbox{mod} \; 3;  \\ \\
\left \lfloor \frac{p}{3} \right \rfloor +3 & \mbox{otherwise}. 
\end{array} \right.
$$  
Therefore, we may assume that $n$ is odd and composite.  The case of $n=15$ can be computed individually, so we also assume that $n \geq 21$.

We observe first that for $$m=\left \lfloor \frac{n}{3} \right \rfloor +2$$ we have 
$$\rho \hat{\;} (\mathbb{Z}_n,m,3) \leq u \hat{\;} (n,m,3) \leq 3m-8 \leq n-2,$$ so we always have
$$\chi \hat{\;} (\mathbb{Z}_n,3) \geq \left \lfloor \frac{n}{3} \right \rfloor +3.$$  

Assume now that $n$ has no prime divisors congruent to 2 mod 3 and that $n$ is divisible by 9; let $m=n/3+3$.  Then $m-1$ and $n$ are relatively prime, since if $d$ is a divisor of both $m-1$ and $n$, then $d$ will divide both $3m-3$ and $n$, and hence also their difference, which is 6.  However, $n$ is odd and $m-1$ is not divisible by 3 (since $m$ is), so $d=1$.  According to Theorem \ref{cor h=3},  
$$\rho \hat{\;} (\mathbb{Z}_n,m,3) \leq \min \{u(n,m,3),3m-10\} \leq 3m-10=n-1,$$ so 
$\chi \hat{\;} (\mathbb{Z}_n,3) \geq   n/3  +4.$

Suppose now that $n$ has a prime divisors congruent to 2 mod 3, and let $p$ be the smallest of these.  We then have
$$\chi \hat{\;} (\mathbb{Z}_n,3) \geq \chi  (\mathbb{Z}_n,3) =v_1(n,3)+1= \left( 1+\frac{1}{p} \right) \frac{n}{3}+1.$$  Now if $n=3p$, then we further have
$$\chi \hat{\;} (\mathbb{Z}_n,3) \geq \left( 1+\frac{1}{p} \right) \frac{n}{3}+2,$$
since for $$m=\left( 1+\frac{1}{p} \right) \frac{n}{3}+1=p+2$$ we have 
$$\rho \hat{\;} (\mathbb{Z}_n,m,3) \leq u \hat{\;} (n,m,3) \leq 3m-8=3p-2=n-2.$$ 
Our proof is now complete.  $\Box$

In \cite{Baj:2013a} we made the following conjecture:

\begin{conj} \label{conj rhohatforh=3}
For all $n$ and $m$, we have equality in Theorem \ref{cor h=3}.

\end{conj}

Correspondingly, we believe that:

\begin{conj} \label{conj for chi hat}
For all values of $n \geq 11$, equality holds in Proposition \ref{h=3 rest crit}.

\end{conj}

We have verified that Conjecture \ref{conj for chi hat} holds for all values of $n \leq 50$, and by Theorems \ref{prime rest crit} and \ref{n even h critical}, it holds when $n$ is prime or even.  As additional support, we prove the following:

\begin{thm}  \label{conj rhohatforh=3 implies conj for chi hat}
Conjecture \ref{conj rhohatforh=3} implies Conjecture \ref{conj for chi hat}.

\end{thm}

{\em Proof}: As we noted before, we may assume that $n$ is odd, composite, and greater than 15.

Suppose first that $n$ has a prime divisor that is congruent to 2 mod 3, and let $p$ be the smallest such prime; since $n$ is odd, $p \geq 5$.  Let us set
$$m=\left( 1+\frac{1}{p} \right) \frac{n}{3}+1.$$
We need to prove that Conjecture \ref{conj rhohatforh=3} implies both of the following statements:

A: $\rho \hat{\;} (\mathbb{Z}_n,m+1,3) =n.$

B: If $\rho \hat{\;} (\mathbb{Z}_n,m,3) <n$, then $n=3p$.   

First, note that $m=\chi (\mathbb{Z}_n, 3)$, so $u(n,m,3)=n$ and thus $u(n,m+1,3)=n$ as well.  Thus, looking at the conjectured formula for $\rho \hat{\;} (\mathbb{Z}_n,m,3)$, to prove statement A, it suffices to verify that 

\noindent A.1: $3(m+1)-3-\gcd(n,(m+1)-1) \geq n$;

\noindent A.2: $3(m+1)-9 \geq n$; and
 
\noindent A.3: If $3(m+1)-10 < n$, then $\gcd(n,(m+1)-1) \neq 7$, $m+1$ is not divisible by 3, and $(m+1)-3$ is not divisible by 5. 

Observe that if $d$ divides both $n$ and $m$, then $d$ divides $3m-n$ as well, and so $$\gcd (n,m) \leq 3m-n=n/p+3,$$ which implies that
$$3(m+1)-3-\gcd(n,(m+1)-1) \geq (p+1) \cdot n/p +3 - (n/p+3)=   n,$$ proving A.1.

To prove A.2, observe that, since $n$ is neither prime nor even, we have $n \geq 3p$, and so
$$3(m+1)-9 =(p+1) \cdot n/p -3 \geq n.$$ 

Similarly, we see that $3(m+1)-10 < n$ may only occur if $n=3p$, in which case $m=p+2$, but then neither $3$ nor $p$ divides $m$, so $\gcd(n,m) =1$; $m+1=p+3$ is not divisible by 3; furthermore, $m-2=p$ is not divisible by 5 (since $p=5$ would give $n=15$, which we excluded).  This proves A.3.

To prove statement B, we will suppose, indirectly, that $n \neq 3p$.  But we assumed that $n$ was odd and composite, so $n=5p$ or $n \geq 7p$; furthermore, if $n=5p$ then, for $p$ to be the smallest prime divisor of $n$ that is congruent to 2 mod 3, $p$ would need to be 5.   For $n=25$ we get $m=11$, but Conjecture \ref{conj rhohatforh=3} implies that $\rho \hat{\;} (\mathbb{Z}_{25},11,3) =25$, so we can rule out $n=25$ and so assume that $n \geq 7p$.   Thus, looking again at the conjectured formula for $\rho \hat{\;} (\mathbb{Z}_n,m,3)$, to prove statement B, it suffices to verify that 

\noindent B.1: $3m-3-\gcd(n,m-1) \geq n$; and

\noindent B.2: If $n \geq 7p$, then $3m-10 \geq n$.

The proofs of B.1 and B.2 are similar to that of A.1 and A.2, respectively---we omit the details.  This completes the proof of statement B.

Assume now that $n$ has no prime divisors congruent to 2 mod 3.  This, of course, means that $n$ itself is not congruent to 2 mod 3. We set $$m=\left \lfloor \frac{n}{3} \right \rfloor +3.$$
We need to prove that Conjecture \ref{conj rhohatforh=3} implies both of the following statements:

C: $\rho \hat{\;} (\mathbb{Z}_n,m+1,3) =n.$

D: If $\rho \hat{\;} (\mathbb{Z}_n,m,3) <n$, then $n$ is divisible by 9.   

This time we have $m=\chi (\mathbb{Z}_n, 3)+2$, so $u(n,m,3)=n$ and thus $u(n,m+1,3)=n$ as well.   Thus, looking at the conjectured formula for $\rho \hat{\;} (\mathbb{Z}_n,m,3)$, to prove statement C, it suffices to verify that 

\noindent C.1: $3(m+1)-3-\gcd(n,(m+1)-1) \geq n$;

\noindent C.2: $3(m+1)-10 \geq n$. 

Suppose that $d$ divides both $n$ and $m$, then $d$ divides $$3m-n=
\left \{
\begin{array}{ll}
9 & \mbox{if } n \equiv 0 \; \mbox{mod} \; 3; \\
8 & \mbox{if } n \equiv 1 \; \mbox{mod} \; 3.
\end{array}
\right.
$$ 
Therefore,
$$3(m+1)-3-\gcd(n,(m+1)-1) \geq 
\left \{
\begin{array}{cl}
n+12-3-9  &   \mbox{if } n \equiv 0 \; \mbox{mod} \; 3; \\
n-1+12-3-8  &   \mbox{if } n \equiv 1 \; \mbox{mod} \; 3.
\end{array}
\right.
$$ This proves C.1.  Since $$m+1 \geq (n-1)/3+ 4,$$  statement C.2 follows as well.

To prove statement D, we first prove that $\gcd(n,m-1) \leq 5$.  Indeed, if $d$ is a divisor of both $n$ and $m-1$, then $d$ divides $3m-3-n$, which is at most 6; however $d$ cannot be 6 as $n$ is odd.  We also see that $$3m-8 \geq n-1+9-8=n.$$  Furthermore, $m \neq 6$ since $n >15$.   

Therefore, according to Conjecture \ref{conj rhohatforh=3},  for  $\rho \hat{\;} (\mathbb{Z}_n,m,3)$ to be less than $n$, we must have either $n$ and $m$ both divisible by 3, or $n$ divisible by $3m-9$ and $m-3$ divisible by 5.  Since in both these cases $n$ is divisible by 3, we have $m=n/3+3$.  We can rule out the second possibility: if $m-3=n/3$ were to be divisible by 5, then $n$ would be as well, contradicting our assumption that $n$ has no prime divisors congruent to 2 mod 3.  This leaves only one possibility: that $n$ and $m$ are both divisible by 3, which implies that $n$ is divisible by 9, as claimed.  Our proof of statement D and thus of Theorem \ref{conj rhohatforh=3 implies conj for chi hat} is now complete.  $\Box$

It is worth mentioning that, as a special case of Conjecture \ref{conj for chi hat}, for odd integers $n \geq 31$, $$\chi \hat{\;} (\mathbb{Z}_n,3) \leq \tfrac{2}{5} n+1.$$
(The additive constant could be adjusted to include odd integers less than 31.)  This conjecture was made by Gallardo, Grekos, et al.~in ~\cite{GalGre:2002a}, and (for large $n$) proved by Lev via the following more general result:

\begin{thm} [Lev; cf.~\cite{Lev:2002a}] \label{Lev on 3-fold}
Let $G$ be an abelian group of order $n$ with $$n \geq 312|L|+923,$$  where, as before, $L$ is the collection of elements of $G$ that have order at most 2.  Then for any subset $A$ of $G$, at least one of the following possibilities holds:
\begin{itemize}
  \item $|A| \leq \tfrac{5}{13}n$; 
  \item $A$ is contained in a coset of an index-two subgroup of $G$;
  \item $A$ is contained in a union of two cosets of an index-five subgroup of $G$; or
  \item $3 \hat{\;} A=G$.
\end{itemize}

\end{thm}

So, in particular, if $n$ is odd, is at least 1235, and a subset $A$ of $\mathbb{Z}_n$ has size more than $2n/5$, then the last possibility must hold, so we get:

\begin{cor} \label{Lev on 3-fold cor}
If $n \geq 1235$ is an odd integer, then $$\chi \hat{\;} (\mathbb{Z}_n,3) \leq \tfrac{2}{5} n+1.$$
\end{cor}
The bound on $n$ in Corollary \ref{Lev on 3-fold cor} can hopefully be reduced.

As another special case of Conjecture \ref{conj for chi hat}, we claim that if $n \geq 83$ is odd and not divisible by five, then $$\chi \hat{\;} (\mathbb{Z}_n,3) \leq \tfrac{4}{11} n+1.$$ 
Theorem \ref{Lev on 3-fold} does not quite yield this: while a careful read of \cite{Lev:2002a} enables us to reduce the coefficient $5/13$ to $(3-\sqrt{5})/2$ (at least for large enough $n$), this is still higher than $4/11$.

It is also worth pointing out that combining Theorem \ref{h crit numb} with Conjecture \ref{conj for chi hat} yields that, when $n \geq 11$, we have
$$\chi  (\mathbb{Z}_n,3)  \leq  \chi \hat{\;} (\mathbb{Z}_n,3) \leq \chi  (\mathbb{Z}_n,3) + 3.$$
 This is in contrast to the fact that for every positive integer $C$, there are values of $n$ and $m$ so that the quantities $\rho \hat{\;} (\mathbb{Z}_n,m,3)$ and $\rho (\mathbb{Z}_n,m,3)$ are further than $C$ away from one another (cf.~\cite{Baj:2013a}).

{\bf Acknowledgments:} The author acknowledges J. Butterworth's and K. Campbell's preliminary work on Theorems \ref{nu} and \ref{h crit numb}, respectively.

\end{document}